\documentclass{article}%
\usepackage{amsfonts}
\usepackage{amssymb}
\usepackage{amsmath}
\usepackage{graphicx}%
\providecommand{\U}[1]{\protect\rule{.1in}{.1in}}

\newtheorem{theorem}{Theorem}

\newtheorem{example}[theorem]{Example}

\begin{document}

\title{Negation and Implication in Partition Logic}
\author{David Ellerman\\University of Ljubljana, Slovenia}
\date{}
\maketitle

\begin{abstract}
\noindent The Boolean logic of subsets, usually presented as `propositional
logic,' is considered as being "classical" while intuitionistic logic and the
many sublogics and off-shoots are "non-classical." But there is another
mathematical logic, the logic of partitions, that is at the same mathmatical
level as Boolean subset logic since subsets and quotient sets (partitions or
equivalence relations) are dual to one another in the category-theoretic
sense. Our purpose here is to explore the notions of negation and implication
in that other mathematical logic of partitions.

\end{abstract}
\tableofcontents

\section{Introduction: The Boolean logic of subsets and the logic of
partitions}

Today, the `classical' form of logic is seen as the Boolean logic of subsets
usually presented as the special case of propositional logic (i.e., the logic
of subsets $0$ and $1$ of the one element set $1$). Other related logics, such
as intuitionistic logic (e.g., the logic of the open subsets of a topological
space) are considered as non-classical. But there is another recently
developed logic that is at the same mathematical level of fundamentality as
subset logic and is thus `classical' in that nontemporal sense. Since the
development of category theory starting in the middle of the twentieth
century, it has been known that the concept of a subset has a
category-theoretic dual in the notion of a quotient set (or, equivalently, a
partition or equivalence relation). Hence, it should be no surprise that there
is a logic of partitions (\cite{ell:lop}; \cite{ell:intropartlogic}) dual to
the Boolean logic of subsets. And since subsets and quotient sets are at the
same basic level from the mathematical point of view, partition logic is more
of a dual sibling to subset logic rather than being another `non-classical'
off-shoot of the classical subset logic.

\section{The logic of partitions}

Our purpose here is briefly present the basics of partition logic that suffice
to explore the role of negation and implication in that logic. A
\textit{partition} $\pi=\left\{  B,B^{\prime},...\right\}  $ on a set $U$ is a
set of non-empty subsets $B$, $B^{\prime}$,... ("blocks") of $U$ where the
blocks are mutually exclusive (the intersection of distinct blocks is empty)
and jointly exhaustive (the union of the blocks is $U$). An
\textit{equivalence relation} is a binary relation $E\subseteq U\times U$ that
is reflexive, symmetric, and transitive. Every equivalence relation on a set
$U$ determines a partition on $U$ where the equivalence classes are the
mutually exclusive and jointly exhaustive blocks of the partition. Conversely,
every partition on a set determines an equivalence relation on the set; two
elements are equivalent if they are in the same block of the partition. The
notions of a partition on a set and an equivalence relation on a set are thus
interdefinable. Indeed, equivalence relations and partitions are often
considered as the "same" as in the conventional practice (not used here) of
defining the "lattice of partitions" as the lattice of equivalence relations
\cite{Birk:lt}.

For the purposes of partition logic, it is important to consider the
complementary binary relation to an equivalence relation. A \textit{partition
relation }$R\subseteq U\times U$\textit{\ }is irreflexive (i.e., $\left(
u,u\right)  \not \in R$ for any $u\in U$), symmetric (i.e., $\left(
u,u^{\prime}\right)  \in R$ implies $\left(  u^{\prime},u\right)  \in R$), and
\textit{anti-transitive}\ in the sense that if $\left(  u,u^{\prime}\right)
\in R$, then for any $a\in U$, either $\left(  u,a\right)  \in R$ or $\left(
a,u^{\prime}\right)  \in R$ (i.e., $U\times U-R=R^{c}$ is transitive). Thus as
binary relations, equivalence relations and partition relations are
complementary. That is, $E\subseteq U\times U$ is an equivalence relation if
and only if (iff) $E^{c}\subseteq U\times U$ is a partition relation.

A \textit{distinction of a partition} is an ordered pair $\left(  u,u^{\prime
}\right)  $ of elements of $U$ in distinct blocks of the partition. The set of
distinctions (abbreviated "dits") of a partition is the \textit{ditset}

\begin{center}
$\operatorname{dit}\left(  \pi\right)  =\left\{  \left(  u,u^{\prime}\right)
:\exists B,B^{\prime}\in\pi;B\neq B^{\prime};u\in B;u^{\prime}\in B^{\prime
}\right\}  $.
\end{center}

\noindent Similarly an \textit{indistinction} or \textit{indit} of a partition
is an ordered pair of elements in the same block of the partition so:

\begin{center}
$\operatorname{indit}\left(  \pi\right)  =\left\{  \left(  u,u^{\prime
}\right)  :\exists B\in\pi;u,u^{\prime}\in B\right\}  =%
{\displaystyle\bigcup\limits_{B\in\pi}}
B\times B=U\times U-\operatorname{dit}\left(  \pi\right)  $.
\end{center}

\noindent The indit set of a partition is the equivalence relation defined by
the partition, and the ditset of a partition is the complementary partition
relation defined by the partition.

If $\sigma=\left\{  C,C^{\prime},...\right\}  $ is another partition on $U$,
then the partial order of \textit{refinement} is defined by:

\begin{center}
$\sigma\precsim\pi$ (read: $\pi$ refines $\sigma$ or $\sigma$ is refined by
$\pi$) if $\forall B\in\pi,\exists C\in\sigma$ such that $B\subseteq C$.
\end{center}

\noindent Note that if $\sigma\precsim\pi$, then for any $C\in\sigma$, there
is a set of blocks of $\pi$ whose union is $C$. The most refined partition on
$U$ is the \textit{discrete partition} $\mathbf{1}=\left\{  \left\{
u\right\}  \right\}  _{u\in U}$ whose blocks are all singletons. It is the top
or maximal element in the refinement partial order. The least refined
partition is the \textit{indiscrete partition} (nicknamed the `blob')
$\mathbf{0}=\left\{  U\right\}  $ whose only block is $U$ itself. It is bottom
or minimal element in the refinement partial order. The \textit{join} $\pi
\vee\sigma$ (least upper bound) of $\pi$ and $\sigma$ is the partition whose
blocks are the non-empty intersections of the blocks of $\pi$ and $\sigma$:

\begin{center}
$\pi\vee\sigma=\left\{  B\cap C\neq\emptyset:B\in\pi;C\in\sigma\right\}  $.
\end{center}

To define the \textit{meet }$\pi\wedge\sigma$ (greatest lower bound) of $\pi$
and $\sigma$, we define an equivalence relation on $U$ that is generated by
$u\sim u^{\prime}$ if $u$ and $u^{\prime}$ are in the same block of $\pi$ or
$\sigma$. Thus if two blocks of $\pi$ and $\sigma$ overlap (non-empty
intersection) then all the elements of the two blocks are equated and so forth
for any finite sequence of overlapping blocks. Hence a block of the meet
partition, i.e., an equivalence class of that equivalence relation, is a
precise union of blocks of $\pi$ and a union of blocks of $\sigma$, and is the
smallest such union. These definitions of refinement, join, and meet turn the
set $\Pi\left(  U\right)  $ of partitions on $U$ into a lattice. The notion of
refinement between partitions is equivalent to inclusion between their
corresponding ditsets or partition relations, i.e., $\sigma\precsim\pi$ iff
$\operatorname{dit}\left(  \sigma\right)  \subseteq\operatorname{dit}\left(
\pi\right)  $, so the lattice of partitions on $U$ can be represented as the
isomorphic lattice of partition relations on $U\times U$. But it should be
carefully noted that what many textbooks call the "lattice of partitions" is
really the opposite lattice of equivalence relations, e.g., Birkhoff
\cite{Birk:lt} or Gr\"{a}tzer \cite{grat:glt}, where the join and meet are interchanged.

The lattice of partitions (in either presentation) was known and studied in
the nineteenth century by Richard Dedekind and others. But no other operations
on partitions besides join and meet were defined throughout the twentieth century.

\begin{quote}
Equivalence relations are so ubiquitous in everyday life that we often forget
about their proactive existence. Much is still unknown about equivalence
relations. Were this situation remedied, the theory of equivalence relations
could initiate a chain reaction generating new insights and discoveries in
many fields dependent upon it.

This paper springs from a simple acknowledgement: the only operations on the
family of equivalence relations fully studied, understood and deployed are the
binary join $\vee$ and meet $\wedge$ operations. \cite[p. 445]{bmp:eqrel}
\end{quote}

\noindent Hence the development of partition logic depended on defining at
least implication $\sigma\Rightarrow\pi$, and then all the other logical
(i.e., Boolean) operations on partitions, e.g., \cite{ell:graph}.

\section{Implication and negation in partition logic}

There are at least four equivalent ways to define the implication operation
$\sigma\Rightarrow\pi$ on partitions. The most intuitive and useful
set-of-blocks definition will be used here. The \textit{implication partition}
$\sigma\Rightarrow\pi$ is like the partition $\pi$ except that every block
$B\in\pi$ that is contained in some block $C\in\sigma$ is replaced by
singletons of its elements. Such an `atomized' or discretized block $B$ might
be denoted $\mathbf{1}_{B}$ as the local $B$-version of the discrete partition
$\mathbf{1}$. If a block $\ B\in\pi$ is not contained in any block of $\sigma
$, then it remains the same which might be denoted $\mathbf{0}_{B}$ as the
local $B$-version of the indiscrete partition $\mathbf{0}$. Hence the
implication partition $\sigma\Rightarrow\pi$ functions as an indicator or
characteristic function with blocks $\mathbf{1}_{B}$ or $\mathbf{0}_{B}$
according to whether or not $B$ was contained in a block of $\sigma$. With the
implication operation, we could refer to $\Pi\left(  U\right)  $ as the
algebra of partitions on $U$ instead of just the lattice of partitions.

If \textit{all} the blocks of $\pi$ are contained in blocks of $\sigma$, i.e.,
if $\pi$ refines $\sigma$, then $\sigma\Rightarrow\pi=\mathbf{1}$, the top of
$\Pi\left(  U\right)  $. Thus we have: $\sigma\Rightarrow\pi=1 $ iff
$\sigma\precsim\pi$ in partition logic just as we have in Boolean subset logic
for the conditional or implication operation ($S\supset T:=S^{c}\cup T$) on
subsets $S,T\subseteq U$: $S\supset T=U$ iff $S\subseteq T$.

With the implication operation, the (absolute) \textit{negation of} $\sigma$
can be defined as $\lnot\sigma:=\sigma\Rightarrow\mathbf{0}$. But the more
interesting (relative)\textit{\ }$\pi$\textit{-negation of }$\sigma$ is
defined as: $\overset{\pi}{\lnot}\sigma:=\sigma\Rightarrow\pi$, so the $\pi
$-negation of $\sigma$ is just another way of considering the implication
$\sigma\Rightarrow\pi$.

The equivalence relation corresponding to the indiscrete partition
$\mathbf{0}$ is the universal relation $U\times U$. For any two equivalence
relations $E,E^{\prime}\subseteq U\times U$, if $E\cup E^{\prime}=U\times U$,
then $E=U\times U$ or $E^{\prime}=U\times U$. This is essentially the standard
result of graph theory that the complement of any disconnected graph is
connected \cite[p. 30]{wil:gt}. Since the indiscrete partition has no
distinctions, i.e., $\operatorname{dit}\left(  \mathbf{0}\right)  =\emptyset$,
the complementary form of that result is that for any two partitions
$\sigma,\pi$, if $\operatorname{dit}\left(  \sigma\right)  \cap
\operatorname{dit}\left(  \pi\right)  =\emptyset$, then $\operatorname{dit}%
\left(  \sigma\right)  =\emptyset$ or $\operatorname{dit}\left(  \pi\right)
=\emptyset$, i.e., $\sigma=\mathbf{0}$ or $\pi=\mathbf{0}$. An alternative
form of the result is useful to understand the negation $\sigma\Rightarrow
\mathbf{0}$.

\begin{theorem}
[Common-Dits Theorem]Any two non-empty ditsets overlap, i.e., have some dits
in common.
\end{theorem}

\noindent\textbf{Proof}: Let $\pi$ and $\sigma$ be any two partitions on $U$
with non-empty dit sets, i.e., $\pi\not =\mathbf{0}\not =\sigma$. We need to
show that $\operatorname*{dit}\left(  \pi\right)  \cap\operatorname*{dit}%
\left(  \sigma\right)  \not =\emptyset$. Since $\sigma$ is not the blob
$\mathbf{0}$, consider two elements $u$ and $u^{\prime}$ distinguished by
$\sigma$ but identified by $\pi$ [otherwise $\left(  u,u^{\prime}\right)
\in\operatorname*{dit}\left(  \pi\right)  \cap\operatorname*{dit}\left(
\sigma\right)  $ and we are finished]. Since $\pi$ is also not the blob, there
must be a third element $u^{\prime\prime}$ not in the same block of $\pi$ as
$u$ and $u^{\prime}$.%

\begin{center}
\includegraphics[
height=1.7469in,
width=1.6016in
]%
{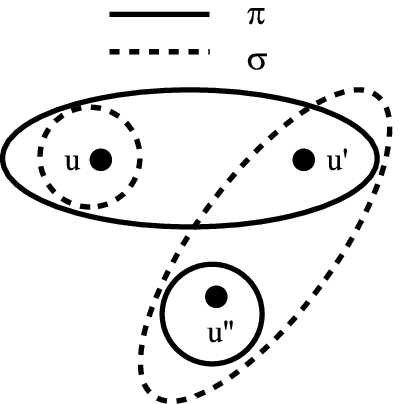}%
\end{center}

\begin{center}
Figure 1: Common dits to any two non-empty ditsets.
\end{center}

\noindent But since $u$ and $u^{\prime}$ are in different blocks of $\sigma$,
the third element $u^{\prime\prime}$ must be distinguished from one or the
other or both in $\sigma$. Hence $\left(  u,u^{\prime\prime}\right)  $ (as in
Figure 1) or $\left(  u^{\prime},u^{\prime\prime}\right)  $ must be
distinguished by both partitions and thus must be in $\operatorname*{dit}%
\left(  \pi\right)  \cap\operatorname*{dit}\left(  \sigma\right)  $. $\square$

This means that for any two non-blob partitions $\pi$ and $\sigma$ on $U$,
there is always a pair of elements $u,u^{\prime}\in U$ that are in different
blocks of both partitions. This result is perhaps particularly striking if we
take $\pi$ and $\sigma$ to be \textit{atomic} partitions, namely, partitions
with only two blocks. For any two ways to divide the elements of $U$
($\left\vert U\right\vert \geq2$) into two parts, there is always a pair of
elements separated by both divisions.

Since intuitionistic logic is the most developed logic aside from Boolean
logic, it is often suggestive to compare the ditsets of partition logic with
the open sets in the topological representation of intuitionistic logic, i.e.,
of a Heyting algebra (also called a pseudo-Boolean algebra or Brouwer
algebra). The negation of an open set is the largest open set disjoint from
the given set. But now we see that there is no non-empty ditsets disjoint from
any given non-empty ditset. Hence intuitively the negation of any partition
$\sigma\neq\mathbf{0}$, is the partition $\mathbf{0}$ with an empty ditset.
The definition $\lnot\sigma:=\sigma\Rightarrow\mathbf{0}$ gives the same
result since the only block $U$ in $\mathbf{0}=\left\{  U\right\}  $ is not
contained in any block of $\sigma\neq\mathbf{0}$. And when $\sigma=\mathbf{0}%
$, then $\lnot\mathbf{0}=\mathbf{0}\Rightarrow\mathbf{0}=\mathbf{1}$ since
$U\subseteq U$ so it is discretized in the implication. That is why the
absolute negation $\lnot\sigma$ is of less interest than the relative negation
$\overset{\pi}{\lnot}\sigma=\sigma\Rightarrow\pi$ which is simply the
partition implication.

\section{Three more equivalent ways to define implication for partitions}

\subsection{The adjunctive definition}

For subsets $R,S,T\subseteq U$, the set implication or conditional $S\supset
T$ in subset logic can also be characterized by a category-theoretic adjunction:

\begin{center}
$R\cap S\subseteq T$ iff $R\subseteq(S\supset T)$.
\end{center}

\noindent Since $S\supset T$ is clearly the maximal subset to satisfy that
characterization, we could define:

\begin{center}
$S\supset T:=\cup\left\{  R:R\cap S\subseteq T\right\}  $.
\end{center}

The partition implication $\sigma\Rightarrow\pi$ can be similarly
characterized in partition logic substituting ditsets for subsets
\cite{ell:lop}. For a third partition $\tau=\left\{  D,D^{\prime},...\right\}
$, the characterization is:

\begin{center}
$\operatorname{dit}\left(  \tau\right)  \cap\operatorname{dit}\left(
\sigma\right)  \subseteq\operatorname{dit}\left(  \pi\right)  $ iff
$\tau\precsim\sigma\Rightarrow\pi$.
\end{center}

Since any intersection of equivalence relations is an equivalence relation,
any union of their complements, the partition relations or ditsets, is also a
ditset. Hence we have a second definition of the partition implication by:

\begin{center}
$\operatorname{dit}\left(  \sigma\Rightarrow\pi\right)  :=\cup\left\{
\operatorname{dit}\left(  \tau\right)  :\operatorname{dit}\left(  \tau\right)
\cap\operatorname{dit}\left(  \sigma\right)  \subseteq\operatorname{dit}%
\left(  \pi\right)  \right\}  $.
\end{center}

\subsection{The graph-theoretic definition}

Another way to define the partition implication or any Boolean operation on
partitions is the graph-theoretic method \cite{ell:graph}. Let $K\left(
U\right)  $ be the complete undirected graph on $U$. The links $u-u^{\prime}$
corresponding to dits of a partition, i.e., $\left(  u,u^{\prime}\right)
\in\operatorname*{dit}\left(  \pi\right)  $, of a partition are labelled with
the `truth value' $T_{\pi}$ and the links corresponding to indits $\left(
u,u^{\prime}\right)  \in\operatorname*{indit}\left(  \pi\right)  $ are
labelled with the `truth value' $F_{\pi}$. Given the two partitions $\pi$ and
$\sigma$, each link in the complete graph $K\left(  U\right)  $ is labelled
with a pair of truth values. Then to define any binary Boolean operation
$\pi\#\sigma$, one evaluates those two truth values on each link of $K\left(
U\right)  $ according to that binary operation to obtain either $T_{\pi
\#\sigma}$ or $F_{\pi\#\sigma}$ on that link. Then we obtain the graph
$G\left(  \pi\#\sigma\right)  $ for that operation by deleting all the links
with the truth value $T_{\pi\#\sigma}$ so that only the $F_{\pi\#\sigma}$
links remain. Those $F_{\pi\#\sigma}$ links then generate an equivalence
relation on $U$ whose blocks are the connected components of the graph
$G\left(  \pi\#\sigma\right)  $. Those connected components or equivalence
classes are the partition $\pi\#\sigma$ on $U$.

Specializing to the implication operation $\sigma\Rightarrow\pi$, the links
retained in $G\left(  \sigma\Rightarrow\pi\right)  $ are the links labelled
with $T_{\sigma}$ and $F_{\pi}$ since that combination is the only one to be
evaluated to $F_{\sigma\Rightarrow\pi}$ in the truth table for implication (or
conditional). The connected components in that graph $G\left(  \sigma
\Rightarrow\pi\right)  $ are the blocks in the partition implication
$\sigma\Rightarrow\pi$.

\begin{example}
Let $U=\left\{  a,b,c,d\right\}  $ so that $K(U)=K_{4}$ is the complete graph
on four points. Let $\sigma=\left\{  \left\{  a\right\}  ,\left\{
b,c,d\right\}  \right\}  $ and $\pi=\left\{  \left\{  a,b\right\}  ,\left\{
c,d\right\}  \right\}  $ so we see immediately from the usual definition, that
the $\pi$-block of $\left\{  c,d\right\}  $ will be discretized while the
$\pi$-block of $\left\{  a,b\right\}  $ will remain whole so the partition
implication is $\sigma\Rightarrow\pi=\left\{  \left\{  a,b\right\}  ,\left\{
c\right\}  ,\left\{  d\right\}  \right\}  $. After labelling the links in
$K\left(  U\right)  $, we see that only the $a-b$ link has the $F_{\sigma
\Rightarrow\pi}$ `truth value' so the graph $G\left(  \sigma\Rightarrow
\pi\right)  $ has only that $a-b$ link (thickened in Figure 2). Then the
connected components of $G\left(  \sigma\Rightarrow\pi\right)  $ give the same
partition implication $\sigma\Rightarrow\pi=\left\{  \left\{  a,b\right\}
,\left\{  c\right\}  ,\left\{  d\right\}  \right\}  $.
\end{example}%

\begin{center}
\includegraphics[
height=1.9233in,
width=4.3751in
]%
{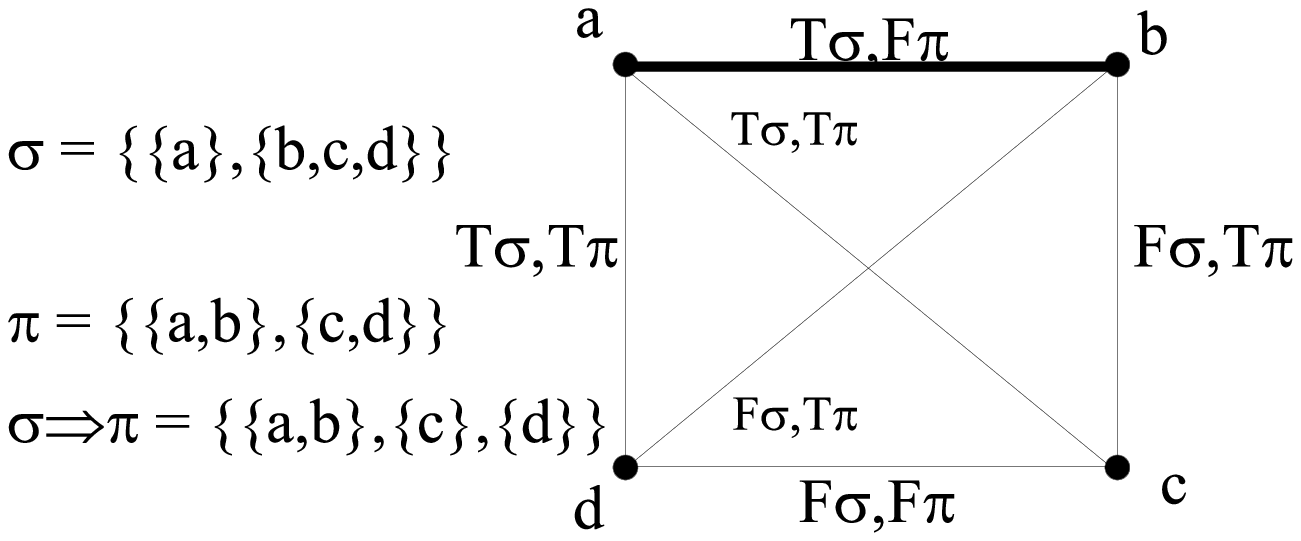}%
\end{center}

\begin{center}
Figure 2: Graph to define the partition implication.
\end{center}

\subsection{The ditset definition}

Another equivalent way to define partition implication is to mimic the subset
definition using ditsets except for the fact that Boolean subset operations on
ditsets do not necessarily lead to ditsets. By the analogy with the interior
operation on subsets of a topological space, the interior $\operatorname*{int}%
\left(  S\right)  $ of an arbitrary subset $S\subseteq U$, where $U$ is a
topological space, is the largest open set contained in $S$. The interior
$\operatorname*{int}\left(  S\right)  $ can also be defined as the complement
of the (topological) closure of the complement, i.e., $\left(  \overline
{S^{c}}\right)  ^{c}$. Similarly, we could start with any subset $S\subseteq
U\times U$, and define the interior of $S$ as the largest ditset contained in
$S$. It could be constructed by first taking the complement $S^{c}$ in
$U\times U$ and then its reflexive, symmetric, and transitive closure
$\overline{S^{c}}$ which is just the intersection of all the equivalence
relations containing $S^{c}$. Then the complement is the \textit{interior}:
$\operatorname*{int}\left(  S\right)  :=\left(  \overline{S^{c}}\right)  ^{c}%
$. This partition-theoretic closure operation is not a topological closure
operation, e.g., since the intersection of two ditsets is not necessarily a
ditset, whereas the intersection of two open sets is open.

We can now use this interior operation to mimic the subset-logic definition of
the implication: $S\supset T:=S^{c}\cup T$. Hence to define the ditset
$\operatorname{dit}\left(  \sigma\Rightarrow\pi\right)  $, we first form
$\operatorname{dit}\left(  \sigma\right)  ^{c}\cup\operatorname{dit}\left(
\pi\right)  $ but that is not a ditset, so the definition is:

\begin{center}
$\operatorname{dit}\left(  \sigma\Rightarrow\pi\right)  =\operatorname*{int}%
\left[  \operatorname{dit}\left(  \sigma\right)  ^{c}\cup\operatorname{dit}%
\left(  \pi\right)  \right]  $.
\end{center}

\noindent Like the graph-theoretic definition, this approach can also be used
for the other operations. The interior operation isn't needed for the ditset
treatment of the join since $\operatorname{dit}\left(  \sigma\vee\pi\right)
=\operatorname{dit}\left(  \sigma\right)  \cup\operatorname{dit}\left(
\pi\right)  $, but the meet could be defined as $\operatorname{dit}\left(
\sigma\wedge\pi\right)  =\operatorname*{int}\left[  \operatorname{dit}\left(
\sigma\right)  \cap\operatorname{dit}\left(  \pi\right)  \right]  $.

Thus we have four definitions of the partition implication $\sigma
\Rightarrow\pi$ that are equivalent (see \cite{ell:lop};
\cite{ell:intropartlogic}).

\section{Relative negation in partition logic}

To study relative negation, we take the `consequence' $\pi$ as fixed and then
let the `antecedent' $\sigma$ vary in the $\pi$-negation $\overset{\pi}{\lnot
}\sigma:=\sigma\Rightarrow\pi$. Another suggestion from intuitionistic logic
is that the negated elements in a Heyting algebra form a Boolean algebra. In
partition logic, this is trivially true for absolute negation since the
negated elements form the two-element Boolean algebra. And it is also true for
the relative $\pi$-negation as was suggested by viewing the implication
$\sigma\Rightarrow\pi$ as an indicator or characteristic function for the
inclusion of the blocks of $\pi$ in the blocks of $\sigma$. And the double
$\pi$-negation $\overset{\pi}{\lnot}\overset{\pi}{\lnot}\sigma=\left(
\sigma\Rightarrow\pi\right)  \Rightarrow\pi$ just interchanges the
$\mathbf{0}_{B}$ and $\mathbf{1}_{B}$ so $\pi$-negation is like the usual
negation of a subset represented by its indicator function (i.e., negation
interchanges the zero-one values). Thus the triple $\pi$-negation is the same
as the single $\pi$-negation. For another partition $\tau=\left\{
D,D^{\prime},...\right\}  $, the join $\overset{\pi}{\lnot}\sigma\vee
\overset{\pi}{\lnot}\tau=\left(  \sigma\Rightarrow\pi\right)  \vee\left(
\tau\Rightarrow\pi\right)  $ would have $B$ discretized, i.e., turned into
$\mathbf{1}_{B}$, iff $B$ is contained in a block $C\in\sigma$ \textit{or} $B$
is contained in a block $D\in\tau$, so it acts like the Boolean join or
disjunction: $\mathbf{0}_{B}\vee\mathbf{1}_{B}=\mathbf{1}_{B}\vee
\mathbf{0}_{B}=\mathbf{1}_{B}\vee\mathbf{1}_{B}=\mathbf{1}_{B}$ and
$\mathbf{0}_{B}\vee\mathbf{0}_{B}=\mathbf{0}_{B}$. Similarly, the meet
$\overset{\pi}{\lnot}\sigma\wedge\overset{\pi}{\lnot}\tau=\left(
\sigma\Rightarrow\pi\right)  \wedge\left(  \tau\Rightarrow\pi\right)  $ would
have $B$ discretized, i.e., turned into $\mathbf{1}_{B}$ iff $B$ is contained
in a block $C\in\sigma$ \textit{and} $B$ is contained in a block $D\in\tau$,
so it acts like the Boolean conjunction. Thus all the partitions over $U$ in
the form of a $\pi$-negation $\overset{\pi}{\lnot}\sigma=\sigma\Rightarrow\pi$
form a Boolean algebra $\mathcal{B}_{\pi}$ with $\pi=\overset{\pi}{\lnot
}\mathbf{1}$ as the bottom element and $\mathbf{1}=\pi\Rightarrow\pi$ as the
top element. Since all the $\pi$-negated partitions (also called $\pi$-regular
partitions) refine $\pi$, the Boolean algebra $\mathcal{B}_{\pi}$ is contained
in the upper segment $\left[  \pi,1\right]  $ and might be called the
\textit{Boolean core} $\mathcal{B}_{\pi}$ of $\left[  \pi,1\right]  $.

There is another construction of $\mathcal{B}_{\pi}$ based on the fact that
singleton blocks in $\pi$ are already atomized so the implication
$\sigma\Rightarrow\pi$ essentially ignores the singletons of $\pi$. Those
singletons are always contained in some block of $\sigma$ so they should be
discretized into singletons, but they are already singletons. If we let
$\pi_{ns}$ stand for the set of non-singleton blocks of $\pi$, then every
$\pi$-negated formula $\sigma\Rightarrow\pi$ is characterized by the set of
non-singleton blocks $B\in\pi_{ns}$ that were discretized, i.e., were assigned
$1_{B}$ instead of $0_{B}$ by the implication $\sigma\Rightarrow\pi$ viewed as
an indicator function (for inclusion of blocks of $\pi$ in blocks of $\sigma
$). Then it is easily seen that the powerset Boolean algebra $\wp\left(
\pi_{ns}\right)  $ on the set of non-singleton blocks of $\pi$ is isomorphic
to the Boolean core $\mathcal{B}_{\pi}$, i.e.,

\begin{center}
$\wp\left(  \pi_{ns}\right)  \cong\mathcal{B}_{\pi}$.
\end{center}

\noindent Thus, we also have: $\wp\left(  \pi\right)  \cong\mathcal{B}_{\pi
}\times%
{\textstyle\prod\nolimits_{\left\{  u\right\}  \in\pi}}
2$.

The single $\pi$-negation $\overset{\pi}{\lnot}\sigma$, the double $\pi
$-negation $\overset{\pi}{\lnot}\overset{\pi}{\lnot}\sigma$, and the excluded
middle $\sigma\vee\overset{\pi}{\lnot}\sigma$ are all partitions of special
interest. The non-singleton blocks of $\overset{\pi}{\lnot}\sigma$ are the
blocks $B\in\pi$ that intersect two or more blocks of $\sigma$. Thus the
non-singleton blocks of the double $\pi$-negation $\overset{\pi}{\lnot
}\overset{\pi}{\lnot}\sigma$ are the blocks $B\in\pi$ that are contained in
blocks of $\sigma$ so $\sigma\Rightarrow\overset{\pi}{\lnot}\overset{\pi
}{\lnot}\sigma$ is a partition tautology and $\sigma\precsim\overset{\pi
}{\lnot}\overset{\pi}{\lnot}\sigma$. The double $\pi$-negation $\overset{\pi
}{\lnot}\overset{\pi}{\lnot}\sigma$ can be thought of as the $\pi
$\textit{-closure} of any $\sigma$ inside of $\mathcal{B}_{\pi}$. Since
$\pi\precsim\overset{\pi}{\lnot}\overset{\pi}{\lnot}\sigma$, we also have that
$\sigma\vee\pi\precsim\overset{\pi}{\lnot}\overset{\pi}{\lnot}\sigma$.

Moreover, since the non-singleton blocks of $\overset{\pi}{\lnot}\sigma$
intersect two or more blocks of $\sigma$, the blocks (always non-singleton
unless otherwise specified) in the excluded middle partition $\sigma
\vee\overset{\pi}{\lnot}\sigma$ are all (strictly) smaller than the blocks of
$\pi$ so $\pi\precsim\sigma\vee\overset{\pi}{\lnot}\sigma$ and thus
$\sigma\vee\pi\precsim\sigma\vee\overset{\pi}{\lnot}\sigma$. And since the
blocks of $\sigma\vee\overset{\pi}{\lnot}\sigma$ are \textit{strictly} smaller
than the blocks of $\pi$, no blocks of $\pi$ are discretized in its $\pi
$-negation, i.e., $\overset{\pi}{\lnot}\left(  \sigma\vee\overset{\pi}{\lnot
}\sigma\right)  =\pi$. Thus the double $\pi$-negation of the excluded-middle
partition is $\mathbf{1}$, i.e., $\overset{\pi}{\lnot}\overset{\pi}{\lnot
}\left(  \sigma\vee\overset{\pi}{\lnot}\sigma\right)  $ is a partition
tautology. While the excluded middle partition $\sigma\vee\overset{\pi}{\lnot
}\sigma$ is not (in general) equal to $\mathbf{1}$ (i.e., is not in general a
partition tautology) nor even in $\mathcal{B}_{\pi}$, it could be said to be
$\pi$\textit{-dense} in $\mathbf{1}$ since its $\pi$-closure is $\mathbf{1}$.

Since both the excluded middle partition $\sigma\vee\overset{\pi}{\lnot}%
\sigma$ and the double $\pi$-negation partition $\overset{\pi}{\lnot}%
\overset{\pi}{\lnot}\sigma$ refine $\sigma\vee\pi$, their meet (greatest lower
bound) $\left(  \sigma\vee\overset{\pi}{\lnot}\sigma\right)  \wedge
\overset{\pi}{\lnot}\overset{\pi}{\lnot}\sigma$ must also refine $\sigma
\vee\pi$. Moreover, that is an equality since the blocks of $\sigma
\vee\overset{\pi}{\lnot}\sigma$ are the non-empty intersections $C\cap B$ for
$C\in\sigma$ and $B\in\pi$ where $B$ is not contained in any $C\in\sigma$, and
the blocks of $\overset{\pi}{\lnot}\overset{\pi}{\lnot}\sigma$ are the blocks
$B$ contained in some $C\in\sigma$. Those non-singleton blocks are all
disjoint, so there are no overlaps in the meet operation. Hence those blocks
remain the same in the meet and they are precisely the blocks of the join
$\sigma\vee\pi$, i.e.,

\begin{center}
$\sigma\vee\pi=\left(  \sigma\vee\overset{\pi}{\lnot}\sigma\right)
\wedge\overset{\pi}{\lnot}\overset{\pi}{\lnot}\sigma$.
\end{center}

\section{Valid formulas}

Since there are partition operations corresponding to all the Boolean
operations \cite{ell:graph}, we can just write logical formulas using those
logical operations without first specifying whether the variables stand for
subsets or partitions (or open subsets as in the intuitionistic case) with the
corresponding operations. Thus we can directly compare the valid formulas in
the different logics.

In the Boolean logic of subsets, a \textit{valid }formula (or \textit{subset
tautology}) is a formula (N.B., a formula, not a proposition) where no matter
what subsets of the universe set $U$ (where $\left\vert U\right\vert \geq1$)
are substituted for the variables, the whole formula evaluates to $U $, the
top of the Boolean algebra of subsets on $U$. The fact that the same set of
valid formulas is obtained if one only considers the two subsets of the
one-element universe $U=1$ is a \textit{theorem} of subset logic (which was
known to Boole). But today most, if not all, textbooks unfortunately ignore
subset logic and present only that special case where $U=1$, "propositional
logic," and then \textit{define} a valid formula as a proposition that is a
truth-table tautology. That common misconception that Boolean logic is just
about propositions (or zero-one entities) rather than subsets seems to have
retarded the development of the dual logic of partitions (since subsets have a
dual, quotient sets or partitions, while propositions do not). Valid formulas
in intuitionistic logic may be defined similarly with open subsets substituted
for arbitrary subsets.

In the logic of partitions on $U$, a \textit{valid} formula or
\textit{partition tautology} is a formula where no matter what partitions on
the universe $U$ (where $\left\vert U\right\vert \geq2$) are substituted for
the variables, the whole formula evaluates to the discrete partition
$\mathbf{1}$, the top of the algebra of partitions on $U$.

There is a simple way to see that all partition tautologies are also subset
tautologies, i.e., valid formulas of subset logic. Consider the partition
algebra $\Pi\left(  2\right)  $ on the two-element set $2=\left\{
0,1\right\}  $. It has only two partitions, the discrete partition
$\mathbf{1}=\left\{  \left\{  0\right\}  ,\left\{  1\right\}  \right\}  $
where $0$ and $1$ are distinguished, and the indiscrete partition
$\mathbf{0}=\left\{  \left\{  0,1\right\}  \right\}  $ where they are not
distinguished. The partition operations, such as join, meet, and implication,
applied to those two partitions could be described in "truth tables" since
$\mathbf{0}$ and $\mathbf{1}$ are the only partitions on $2$. And those truth
tables are the same as the Boolean subset operations on the two subsets of the
one-element set. Hence we have an isomorphism between the partition algebra
$\Pi\left(  2\right)  $ on $2$ and the power-set Boolean algebra $\wp\left(
1\right)  =2$ for $1$ as the one-element set. Now consider any formula that is
a valid formula in partition logic. Since it evaluates to $\mathbf{1}$ for all
partitions on any $U$ where $\left\vert U\right\vert \geq2$, it does that for
$U=2$, but then the isomorphism $\Pi\left(  2\right)  \cong\wp\left(
1\right)  $ means that the same formula will be a truth table tautology in
$\wp\left(  1\right)  $ and thus it is a valid formula for subset logic in
general. Hence all partition tautologies are subset tautologies. But the
inclusion is strict. For instance, the law of excluded middle $\sigma\vee
\lnot\sigma=\sigma\vee\left(  \sigma\Rightarrow\mathbf{0}\right)  $ is not a
partition tautology since for any $\sigma\neq\mathbf{0}$, $\mathbf{1} $,
$\sigma\Rightarrow\mathbf{0}=\mathbf{0}$, and $\sigma\vee\lnot\sigma
=\sigma\vee\mathbf{0}=\sigma\neq\mathbf{1}$.

The Boolean core $\mathcal{B}_{\pi}$ of the upper segment $\left[
\pi,\mathbf{1}\right]  $ for any partition $\pi$, provides a way to
`automatically' generate partition tautologies. Since $\mathcal{B}_{\pi}$ is a
Boolean algebra, any Boolean tautology comprised of $\pi$-negated partitions
will also be a partition tautology. For instance, the law of excluded middle
in $\mathcal{B}_{\pi}$ has the form $\overset{\pi}{\lnot}\sigma\vee
\overset{\pi}{\lnot}\overset{\pi}{\lnot}\sigma$ which is the "weak law of
excluded middle" in partition logic. It is a partition validity since it is a
Boolean tautology that evaluates to $\mathbf{1}$ no matter what partitions on
$U$ are substituted for $\pi$ and $\sigma$.

Conversely, given any formula using the connectives of $\vee$, $\wedge$,
$\Rightarrow$, and the constants of $\mathbf{0}$ and $\mathbf{1}$, its
\textit{single }$\pi$\textit{-negation transform} is obtained by replacing
each atomic variable $\sigma$ by its single $\pi$-negation $\overset{\pi
}{\lnot}\sigma=\sigma\Rightarrow\pi$ and by replacing the constant
$\mathbf{0}$ by $\pi$. The binary operations $\vee$, $\wedge$, and
$\Rightarrow$ as well as the constant $\mathbf{1}$ all remain the same. For
instance, the single $\pi$-negation transform of the excluded middle formula
$\sigma\vee\lnot\sigma=\sigma\vee\left(  \sigma\Rightarrow\mathbf{0}\right)  $
is the partition tautology of the weak law of excluded middle for $\pi$-negation:

\begin{center}
$\left(  \sigma\Rightarrow\pi\right)  \vee\left(  \left(  \sigma\Rightarrow
\pi\right)  \Rightarrow\pi\right)  =\overset{\pi}{\lnot}\sigma\vee\overset
{\pi}{\lnot}\overset{\pi}{\lnot}\sigma$.
\end{center}

\noindent Then the single $\pi$-negation transform of any classical tautology
will still be a tautology but now expressed in $\mathcal{B}_{\pi}$ and thus it
is also a partition tautology. Thus all partition tautologies are ordinary
Boolean logic tautologies, and any ordinary subset tautology transforms into a
partition tautology via the single $\pi$-negation transform.

The weak law of excluded middle is also an example of a partition tautology
that is not an intuitionistic validity. Since the lattice of partitions is the
standard example of a non-distributive lattice while intuitionistic logic or
Heyting algebras are distributive, the distributive laws are examples of
formulas that are valid in intuitionistic logic but not in partition logic.
Thus there is no inclusion either way between partition and intuitionistic tautologies.

Since partition lattices and their non-distributivity were known to Dedekind
and some other European mathematicians such as Ernest Schr\"{o}der, it was an
embarrassing moment in American mathematics when the philosopher-mathematican,
Charles Saunders Peirce, claimed to prove the distributivity of all lattices
\cite{peirce:logic} but omitted the `proof' as being too tedious. Europeans
soon besieged him with examples of partition lattices including the simplest
non-trivial one on a three-element set $U=\left\{  a,b,c\right\}  $.%

\begin{center}
\includegraphics[
height=1.5298in,
width=2.258in
]%
{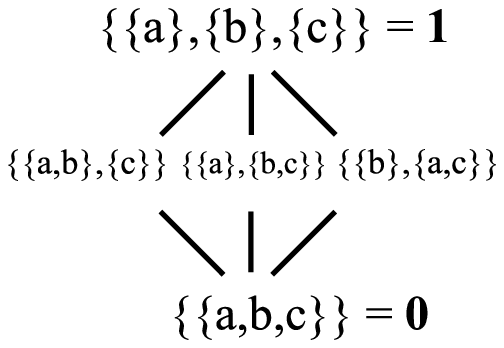}%
\end{center}

\begin{center}
Figure 3: Partition lattice on $U=\left\{  a,b,c\right\}  $.
\end{center}

\noindent Taking the three middle partitions, $\pi=\left\{  \left\{
a,b\right\}  ,\left\{  c\right\}  \right\}  $, $\sigma=\left\{  \left\{
a\right\}  ,\left\{  b,c\right\}  \right\}  $, and $\tau=\left\{  \left\{
b\right\}  ,\left\{  a,c\right\}  \right\}  $, then: $\pi\vee\left(
\sigma\wedge\tau\right)  =\pi\vee\mathbf{0}=\pi$ and $\left(  \pi\vee
\sigma\right)  \wedge\left(  \pi\vee\tau\right)  =\mathbf{1}\wedge
\mathbf{1}=\mathbf{1}$. If Peirce had known about the partition implication
and the Boolean core $\mathcal{B}_{\pi}$, then he could at least have pointed
out that the Boolean core $\mathcal{B}_{\pi}$ for any $\pi$ is distributive,
and moreover any partition $\varphi\in\left[  \pi,\mathbf{1}\right]  $
distributes across the Boolean core \cite{ell:intropartlogic} in the sense that:%

\begin{align*}
\varphi\vee\left(  \overset{\pi}{\lnot}\sigma\wedge\overset{\pi}{\lnot}%
\tau\right)   & =\left(  \varphi\vee\overset{\pi}{\lnot}\sigma\right)
\wedge\left(  \varphi\vee\overset{\pi}{\lnot}\tau\right) \\
\varphi\wedge\left(  \overset{\pi}{\lnot}\sigma\vee\overset{\pi}{\lnot}%
\tau\right)   & =\left(  \varphi\wedge\overset{\pi}{\lnot}\sigma\right)
\vee\left(  \varphi\wedge\overset{\pi}{\lnot}\tau\right)  \text{.}%
\end{align*}

\section{Concluding remarks}

Our purpose has been to develop the notions of negation and implication
(relative negation) in the logic of partitions. Since partition relations
(ditsets) and equivalence relations (indit sets) are complementary in $U\times
U$, every result in the logic of partitions has a complementary-dual result in
the logic of equivalence relations (\cite{ell:lop}; \cite{ell:intropartlogic})
so the latter is not really a different logic but a complementary way to view
partition logic. There is a similar complementary-duality in intuitionistic
logic between Heyting algebras (modelled by the open subsets of a topological
space) and Co-Heyting algebras (modelled by the closed subsets).
Intuitionistic logic makes the symmetry-breaking choice to deal with Heyting
algebras rather than Co-Heyting algebras, and we have made the similar choice
to develop the logic of partitions rather than the (`anti-isomorphic') logic
of equivalence relations. For instance, the complementary-dual to the
implication operation on partitions is the difference operation on equivalence
relations. The partition logic tautology of \textit{modus ponens} has the
customary form: $\left(  \sigma\wedge\left(  \sigma\Rightarrow\pi\right)
\right)  \Rightarrow\pi$, whereas the corresponding formula in the dual logic
of equivalence relations is the unfamiliar $\pi-\left(  \sigma\vee\left(
\pi-\sigma\right)  \right)  $. Hence we have made no independent development
of the equivalence relation notions of difference or of negation as
"difference from $U\times U$."

No new logical operations on partitions, aside from join and meet, were
defined throughout the twentieth century. The definition of the partition
implication (or relative negation) in any of the many equivalent ways was the
key to the development of the full logic of partitions. Why the delay? One
reason is perhaps the fact that partition lattices are so general that any
partition tautologies or identities involving just the lattice operations and
top and bottom, e.g., $\mathbf{1}\wedge\pi=\pi$ or $\mathbf{0}\vee\pi=\pi$,
are in fact identities that hold on all such lattices \cite{whitman:lat}. Thus
the logic of general set partitions only becomes of independent interest by
moving beyond the lattice operations on partitions. Alternatively, one could
develop a `logic' of equivalence relations sticking with only the lattice
operations but specializing to certain types of equivalence relations
\cite{rota:logiccers}. But the main reason for the delay seems to be that the
Boolean logic of subsets was and is presented in only the special case of the
logic of propositions, and propositions, unlike subsets, do not have a
category-theoretic dual concept. Hence twentieth century mathematical
logicians were not even looking for the dual logic of quotient sets,
equivalence relations, or partitions.

\end{document}